\documentclass[a4paper]{amsart}
\usepackage[utf8]{inputenc}
\usepackage[T2A]{fontenc}
\usepackage[english]{babel}

\usepackage{amssymb}

\newtheorem{lemma}{Lemma}
\newtheorem{theorem}{Theorem}

\def\Z{{\mathbb Z}}
\def\C{{\mathbb C}}
\def\Zi{{\mathbb Z}{[i]}}
\def\taue{\tau^{(e)}}
\def\tauec{{\mathfrak t}^{(e)}}
\def\p{{\mathfrak p}}
\def\q{{\mathfrak q}}
\def\dd{{\mathfrak d}}
\def\tt{{\mathfrak t}}
\def\eps{\varepsilon}
\def\res{\mathop{\mathrm{res}}}
\def\llog{\mathop{\mathrm{llog}}}
\def\lllog{\mathop{\mathrm{lllog}}}
\def\const{\mathop{\mathrm{const}}}
\def\suma{\mathop{\sum\nolimits'}}
\def\proda{\mathop{\prod\nolimits'}}
\def\le{\leqslant}
\def\ge{\geqslant}
\def\proof{\par {\em Proof.}\hspace{1em}}
\def\endproof{{\hfil$\blacksquare$\parfillskip0pt\par\medskip}}

\usepackage{hyperref}

\begin{document}

\title[Multidimensional e-divisor function over $\Zi$]{Multidimensional exponential divisor function \\ over Gaussian integers}
\author{Andrew V. Lelechenko}
\address{I.~I.~Mechnikov Odessa National University}
\email{1@dxdy.ru}


\begin{abstract}
Let $\taue_k \colon \Z\to\Z$ be a multiplicative function such that $ \taue_k(p^a) = \sum_{d_1\cdots d_k=a} 1 $. In the present paper we introduce generalizations of $\taue_k$ over the ring of Gaussian integers $\Zi$. We determine their maximal orders by proving a general result and establish asymptotic formulas for their average orders.
\end{abstract}

\maketitle

\section{Introduction}

Exponential divisor function $\taue \colon \Z\to\Z$ introduced by Subbarao in \cite{subbarao1972} is a multiplicative function such that
$$ \taue(p^a) = \tau(a), $$
where $\tau\colon \Z\to\Z$ stands for the usual divisor function, $\tau(n) = \sum_{d|n} 1$. Erdős estimated its maximal order and Subbarao proved an asymptotic formula for $ \sum_{n\le x} \taue(n)$. Later Wu \cite{wu1995} gave more precise estimate:
$$ \sum_{n\le x} \taue(n) = Ax + Bx^{1/2} + O(x^{\theta_{1,2}+\eps}), $$
where $A$ and $B$ are computable constants, $\theta_{1,2}$ is an exponent in the error term of the estimate $\sum_{ab^2 \le x} 1 = \zeta(2) x + \zeta(1/2) x^{1/2} + O(x^{\theta_{1,2}+\eps})$. The best modern result~\cite{graham1988} is that $\theta_{1,2} \le 1057/4785$.

One can consider multidimensional exponential divisor function $\taue_k \colon \Z\to\Z$ such that
$$ \taue_k(p^a) = \tau_k(a), $$
where $\tau_k(n)$ is a number of ordered $k$-tuples of positive integers $(d_1,\ldots,d_k)$ such that $d_1\cdots d_k=n$. So $\taue \equiv \taue_2$. Toth \cite{toth2007a} investigated asymptotic properties of $\taue_k$ and proved that for arbitrarily $\eps>0$
$$ \sum_{n\le x} \taue_k(n) = C_k x + x^{1/2} S_{k-2}(\log x) + O(x^{w_k+\eps}), $$
where $S_{k-2}$ is a polynomial of degree $k-2$ and $w_k=(2k-1)/(4k+1)$.

In the present paper we generalize multidimensional exponential divisor function over the ring of Gaussian integers $\Zi$. Namely we introduce multiplicative functions  $\taue_{k*}\colon\Z\to\Z$, $\tauec_k, \tauec_{k*} \colon \Zi\to\Z$ such that
\begin{equation}\label{eq:definition}
\taue_{k*}(p^a) = \tt_k(a), \qquad
\tauec_k(\p^a) = \tau_k(a), \qquad \tauec_{k*}(\p^a)=\tt_k(a),
\end{equation}
where $p$ is prime over $\Z$, $\p$ is prime over $\Zi$, $\tt_k(a)$ is a number of ordered $k$-tuples of non-associated in pairs Gaussian integers $(\dd_1,\ldots,\dd_k)$ such that $\dd_1\cdots \dd_k=a$.

The aims of this paper are to determine maximal orders of $\taue_k$, $\taue_{k*}$, $\tauec_k$, $\tauec_{k*}$ and to provide asymptotic formulas for $ \sum_{n\le x} \taue_{k*}(n)$, $ \suma_{\!\!\!N(\alpha)\le x} \tauec_k(\alpha)$, $ \suma_{\!\!\!N(\alpha)\le x} \tauec_{k*}(\alpha)$. A theorem on the maximal order of multiplicative functions over $\Zi$, generalizing \cite{suryanarayana1975}, is also proved.

\section{Notation}

Let us denote the ring of Gaussian integers by $\Zi$, $N(a+bi) = a^2+b^2$.

In asymptotic relations we use $\sim$, $\asymp$, Landau symbols $O$ and $o$, Vinogradov symbols $\ll$ and $\gg$ in their usual meanings. All asymptotic relations are written for the argument tending to the infinity.

Letters $\p$ and $\q$ with or without indexes denote Gaussian primes; $p$ and $q$ denote rational primes.

As usual $\zeta(s)$ is Riemann zeta-function and $L(s,\chi)$ is Dirichlet $L$-function. Let~$\chi_4$ be the single nonprincipal character modulo 4, then
$$Z(s) = \zeta(s) L(s, \chi_4)$$
is Hecke zeta-function for the ring of Gaussian integers.

Real and imaginary components of the complex $s$ are denoted as $\sigma:=\Re s$ and~$t:=\Im s$, so $s=\sigma+it$.

We use abbreviations $\llog x := \log\log x$, $\lllog x := \log\log\log x$.

Notation $\suma$ means a summation over non-associated elements of $\Zi$, and $\proda$ means the similar relative to multiplication. Notation $a\sim b$ means that~$a$ and $b$ are associated, that is $a/b \in \{\pm1, \pm i\}$. But in asymptotic relations $\sim$ preserve its usual meaning.

Letter $\gamma$ denotes Euler–Mascheroni constant. Everywhere $\eps>0$ is an arbitrarily small number (not always the same).

We write $f\star g$ for the notation of the Dirichlet convolution
$$ (f \star g)(n) = \sum_{d|n} f(d) g(n/d). $$

\section{Preliminary lemmas}

We need following auxiliary results.

\begin{lemma}[Gauss criterion]\label{l:gaussian-primes}
Gaussian integer $\p$ is prime if and only if one of the following cases complies:
\begin{itemize}
\item $\p \sim 1+i$,
\item $\p \sim p$, where $p \equiv 3 \pmod 4$,
\item $N(\p) = p$, where $p \equiv 1 \pmod 4$.
\end{itemize}
In the last case there are exactly two non-associated $\p_1$ and~$\p_2$ such that $N(\p_1) = N(\p_2) = p$.
\end{lemma}

\proof
See \cite[\S34]{gauss1832}.
\endproof

\begin{lemma}\label{l:prime-rule}
\begin{eqnarray}
\label{eq:prime-rule}
&&\suma_{N(\p)\le x} 1 \sim {x\over\log x},
\\
\label{eq:prime-log-sum}
&&\suma_{N(\p)\le x} \log N(\p) \sim x,
\end{eqnarray}
\end{lemma}

\proof
Taking into account Gauss criterion and the asymptotic law of the distribution of primes in the arithmetic progression we have
\begin{multline*}
\suma_{N(\p)\le x} 1 \sim \#\{ p \mid p\equiv3 \!\!\!\!\!\pmod4, p\le \sqrt{x}  \}
+ 2\#\{ p \mid p\equiv1 \!\!\!\!\!\pmod4, p\le x  \}
\sim\\
\sim {\sqrt x \over \phi(4)\log x/2} + 2 {x\over \phi(4)\log x} = {x\over\log x}.
\end{multline*}
A partial summation with the use of \eqref{eq:prime-rule} gives us the second statement of the lemma.
\endproof

\begin{lemma}\label{l:max-log-tau}
For $k\ge2$
\begin{equation}
\label{eq:max-log-tau}
\max_{n\ge1} {\log \tau_k(n)\over n} = {\log k\over 2},
\end{equation}
\end{lemma}

\proof
Taking into account
$$ \tau_k(p^a) = { k+a-1 \choose a } \le k^a, $$
for $\Omega(n):= \sum_{p^a||n} a$ we have
$\tau_k(n) \le k^{\Omega(n)} \le k^{\log_2 n}$. This implies
$$
{\log \tau_k(n) \over n} \le
{\log_2 n\over n} \log k \le
{\log k \over 2},
$$
because $n^{-1} \log_2 n$ is strictly decreasing for $n\ge2$. But
$$ {\log \tau_k(2)\over2} = {\log k \over 2}. $$
\endproof

\begin{lemma}\label{l:max-log-tau*}
For $k\ge2$
\begin{equation}
\label{eq:max-log-tau*}
\max_{n\ge1} {\log \tt_k(n) \over n} = {1\over2} {\log {k+1 \choose 2}}.
\end{equation}
\end{lemma}

\proof
Let $k_2:= {k+1 \choose 2}$. Lemma \ref{l:gaussian-primes} implies that
\begin{eqnarray*}
&\tt_k(2^a) = \displaystyle{k+2a-1 \choose 2a} \le k_2^a, &\\
&\tt_k(p^a) = \displaystyle{k+a-1 \choose a} \le k^a \le k_2^a &\quad \text{if~} p\equiv 3\pmod 4, \\
&\tt_k(p^a) = \displaystyle{k+a-1 \choose a}^2 \le k^{2a} &\quad \text{if~} p\equiv 1\pmod 4.
\end{eqnarray*}
Let us define
$$ \Omega_1(n):= \sum_{p^a||n \atop p\equiv 1\!\!\!\!\!\pmod 4} a, \quad \Omega_2(n):= \sum_{p^a||n \atop p\not\equiv 1\!\!\!\!\!\pmod 4} a. $$
Then
$ \tt_k(n) \le k^{2\Omega_1(n)} k_2^{\Omega_2(n)} $.
Consider
$$ f(x,y) = {x \log k^2 + y \log k_2 \over 5^x 2^y}, $$
then
$ {n^{-1} \log \tt_k(n)} \le f(\Omega_1(n), \Omega_2(n)) $. One can verify that if $x\ge 1$ or $y\ge 1$ then
$$
f(x+1, y) \le f(x,y),
\qquad
f(x, y+1) \le f(x,y),
$$
because $ \log k_2 + \log k^2 < 5 \log k_2 $. So
$$ \max_{x,y\ge0} f(x,y) = \max \bigl\{ f(1,0), f(0,1) \bigr\} = {\log k_2 \over 2}.$$
But $$ {\log \tt_k(2)\over2} = {\log k_2 \over 2}. $$
\endproof

\begin{lemma}\label{l:rational-maximal-order}
Let $F\colon \Z\to\C$ be a multiplicative function such that  $F_k(p^a)=f(a)$, where $f(n) \ll n^\beta$ for some $\beta>0$. Then
\begin{equation*}\label{eq:rational-maximal-order}
\limsup_{n\to\infty} {\log F_k(n)\llog n \over\log n} = \sup_{n\ge1} {\log f(n)\over n}.
\end{equation*}
\end{lemma}

\proof
See \cite{suryanarayana1975}.
\endproof

\begin{lemma}\label{l:log-sum}
Let $f(t)\ge 0$. If
$$ \int_1^T f(t) \, dt \ll g(T), $$
where $g(T) = T^\alpha \log^\beta T$, $\alpha\ge 1$,
then
\begin{equation*}\label{eq:log-summing}
I(T):= \int_1^T {f(t)\over t} dt \ll
\left\{ \begin{matrix}
\log^{\beta+1} T & \text{if } \alpha=1, \\
T^{\alpha-1} \log^{\beta} T & \text{if } \alpha>1.
\end{matrix} \right.
\end{equation*}
\end{lemma}

\proof
Let us divide the interval of integration into parts:
$$
I(T)
\le
\sum_{k=0}^{\log_2 T}
\int_{T/2^{k+1}}^{T/2^k} {f(t)\over t} dt
<
\sum_{k=0}^{\log_2 T} {1\over T/2^{k+1}}
\int_1^{T/2^k} f(t) dt
\ll
\sum_{k=0}^{\log_2 T} {g(T/2^{k})\over T/2^{k+1}}.
$$
Now the lemma's statement follows from elementary estimates.
\endproof

\begin{lemma}\label{l:zeta-moments}
Let $T>10$ and $|d- 1/2| \ll 1/\log T$. Then we have the following estimates
\begin{eqnarray*}
\int_{d-iT}^{d+iT} |\zeta(s)|^4 {ds\over s} &\ll& \log^5 T,
\\
\int_{d-iT}^{d+iT} |L(s,\chi_4)|^4 {ds\over s} &\ll& \log^5 T,
\end{eqnarray*}
for growing $T$.
\end{lemma}

\proof
The statement is the result of the application of Lemma \ref{l:log-sum} to the estimates \cite[Th. 10.1, p. 75]{montgomery1971}.

\begin{lemma}\label{l:phragmen}
Define $\theta>0$ such that $\zeta(1/2+it) \ll t^\theta$ as $t\to\infty$, and let~$\eta>0$ be arbitrarily small. Then
$$
\zeta(s) \ll \left\{ \begin{matrix}
|t|^{1/2 - (1-2\theta)\sigma}, & \sigma\in[0, 1/2],
\\
|t|^{2\theta(1-\sigma)} , & \sigma\in[1/2, 1-\eta], \\
|t|^{2\theta(1-\sigma)} \log^{2/3} |t| , & \sigma\in[1-\eta, 1], \\
\log^{2/3} |t|, & \sigma\ge1.
\end{matrix}\right.
$$
The same estimates are valid for $L(s, \chi_4)$ also.
\end{lemma}

\proof
The statement follows from Phragmén---Lindelöf principle, exact and approximate functional equations for $\zeta(s)$ and $L(s,\chi_4)$. See \cite{ivic2003} and \cite{titchmarsh1986} for details.
\endproof

The best modern result \cite{huxley2005} is that $\theta \le 32/205+\eps$.

\section{Main results}

First we give maximal orders of $\taue_k$, $\taue_{k*}$, $\tauec_k$ and $\tauec_{k*}$.

The following theorem generalizes Lemma \ref{l:rational-maximal-order} to Gaussian integers; the proof's outline follows the proof of Lemma \ref{l:rational-maximal-order} in \cite{suryanarayana1975}.

\begin{theorem}\label{th:gaussian-maximal-order}
Let $F\colon \Zi\to\C$ be a multiplicative function such that $F(\p^a)=f(a)$, where $f(n) \ll n^\beta$ for some $\beta>0$. Then
\begin{equation*}\label{eq:complex-maximal-order}
\limsup_{\alpha\to\infty} {\log F(\alpha)\llog N(\alpha) \over\log N(\alpha)} = \sup_{n\ge1} {\log f(n)\over n} := K_f.
\end{equation*}
\end{theorem}

\proof
Let us fix arbitrarily small $\eps>0$.

Firstly, let us show that there are infinitely many $\alpha$ such that
$$ {\log F(\alpha) \llog N(\alpha) \over \log N(\alpha)} > K_f - \eps. $$
By definition of $K_f$ we can choose $l$ such that
$$\bigl( \log f(l) \bigr) / l > K_f - \eps/2.$$
It follows from \eqref{eq:prime-log-sum} that for $x\ge2$ inequality
$$\suma_{\!\!\!N(\p)\le x} \log N(\p) > Ax$$
holds, where~$0<A<1$.

Let $\q$ be an arbitrarily large Gaussian prime, $N(\q)\ge2$. Consider
$$
r = \suma_{N(\p)\le N(\q)} 1,
\qquad
\alpha = \proda_{N(\p)\le N(\q)} \p^l.
$$

Then $F(\alpha) = \bigl( f(l) \bigr)^r$ and we have
\begin{equation}\label{eq:max-order-proof*}
r \log N(\q) \ge {\log N(\alpha) \over l}
= \suma_{N(\p)\le N(\q)} \log N(\p) > A N(\q),
\end{equation}
\begin{equation}\label{eq:max-order-proof**}
\log F(\alpha) = r \log f(l) \ge {\log N(\alpha) \over \log N(\q)} {\log f(l) \over l}.
\end{equation}
But  \eqref{eq:max-order-proof*} implies
$$ \log A + \log N(\q) < \log {\log N(\alpha) \over l} \le \llog N(\alpha), $$
so $\log N(\q) < \llog N(\alpha) - \log A$.
Then it follows from \eqref{eq:max-order-proof**} that
$$ \log F(\alpha) > {\log N(\alpha) \over \llog N(\alpha) - \log A} {\log f(l) \over l}  $$
and since $\bigl( \log f(l) \bigr) / l > K_f - \eps/2$ and $A<1$ we have
$$ {\log F(\alpha) \llog N(\alpha) \over \log N(\alpha)} > {\llog N(\alpha) \over \llog N(\alpha) - \log A} (K_f - \eps/2) > K_f - \eps. $$

\medskip

Second, let us show the existence of $N(\eps)$ such that for all $n\ge N(\eps)$ we have
$$
{\log F(n) \llog N(\alpha) \over \log N(\alpha)} < (1+\eps) K_f.
$$
Let us choose $\delta\in(0,\eps)$ and $\eta\in \bigl(0, \delta/(1+\delta) \bigr)$. Suppose $N(\alpha)\ge3$, then we define
$$
\omega := \omega(\alpha) = {(1+\delta) K_f \over \llog N(\alpha)},
\qquad
\Omega := \Omega(\alpha) = \log^{1-\eta} N(\alpha).
$$
By choice of $\delta$ and $\eta$ we have
$$ \Omega^\omega = \exp(\omega\log\Omega) = \exp\bigl( (1-\eta)(1+\delta) K_f \bigr) > e^{K_f}.$$
Suppose that the canonical expansion of $\alpha$ is $$\alpha \sim \p_1^{a_1}\cdots p_r^{a_r} \q_1^{b_1} \cdots \q_s^{b_s},$$ where $N(\p_k) \le \Omega$ and  $N(\q_k) > \Omega$. Then
\begin{equation}\label{eq:max-order-proof***}
{F(\alpha) \over N^\omega(\alpha)}
=
\prod_{k=1}^r {f(a_k) \over N^{\omega a_k} (\p_k)}
\cdot
\prod_{k=1}^s {f(b_k) \over N^{\omega b_k} (\q_k)}
:= \Pi_1 \cdot \Pi_2
\end{equation}
But since $\Omega^\omega > e^{K_f}$  and  $K_f \ge \bigl( \log f(b_k) \bigr) / b_k$ then
$$
{f(b_k) \over N^{\omega b_k}(q_k)}
<
{f(b_k) \over \Omega^{\omega b_k}}
<
{f(b_k) \over e^{K_f b_k}}
\le 1
$$
and it follows that $\Pi_2 \le 1$. Consider $\Pi_1$. From the statement of the theorem we have $f(n) \ll n^\beta$, so
$$
{f(a_k) \over N^{\omega a_k} (p_k)}
\ll
{a_k^\beta \over (\omega a_k )^{\beta} }
\ll
\omega^{-\beta}.
$$
Then
$$ \log \Pi_1 \ll \Omega \log w^{-\beta}
\ll \log^{1-\eta} N(\alpha) \lllog N(\alpha)
= o\left( \log N(\alpha) \over \llog N(\alpha) \right)
$$
And finally by \eqref{eq:max-order-proof***} we get
$$
\log F(n) = \omega \log n + \log \Pi_1 + \log \Pi_2
= {(1+\delta) K_f \log n \over \llog n} + {(\eps-\delta) K_f \log n \over \llog n}.
$$
\endproof

\begin{theorem}\label{th:max-order}
\begin{eqnarray*}
\label{eq:max-order-real}
\limsup_{n\to\infty}
{\log \taue_k(n) \llog n \over \log n} &=& {\log k\over 2},
\\
\label{eq:max-order*-real}
\limsup_{n\to\infty}
{\log \taue_{k*}(n) \llog n \over \log n} &=& {1\over2} {\log {k+1 \choose 2}},
\\
\label{eq:max-order}
\limsup_{\alpha\to\infty}
{\log \tauec_k(\alpha) \llog N(\alpha) \over \log N(\alpha)} &=& {\log k \over2},
\\
\label{eq:max-order*}
\limsup_{\alpha\to\infty}
{\log \tauec_{k*}(\alpha) \llog N(\alpha) \over \log N(\alpha)} &=& {1\over2} {\log {k+1 \choose 2}}.
\end{eqnarray*}
\end{theorem}

\proof
Statements follow from \eqref{eq:max-log-tau}, \eqref{eq:max-log-tau*}, Lemma \ref{l:rational-maximal-order} and Theorem \ref{th:gaussian-maximal-order}.
\endproof

A simple corollary of the Theorem \ref{th:max-order} is that
\begin{equation}\label{eq:max-order-simple}
\taue_{k*}(n) \ll n^\eps,
\qquad
\tauec_k(\alpha) \ll N^\eps(\alpha),
\qquad
\tauec_{k*}(\alpha) \ll N^\eps(\alpha).
\end{equation}

\bigskip

We are ready to provide asymptotic formulas for sums of $ \taue_{k*}(n)$, $  \tauec_k(\alpha)$, $ \tauec_{k*}(\alpha)$.

Let us denote
\begin{eqnarray*}
G_{k*}(s):= \sum_n \taue_{k*}(n) n^{-s},
&\qquad&
T_{k*}(x):= \sum_{n \le x} \taue_{k*}(n),
\\
F_k(s):= \suma_\alpha \tauec_k(\alpha) N^{-s}(\alpha),
&\qquad&
M_k(x):= \suma_{N(\alpha) \le x} \tauec_k(\alpha),
\\
F_{k*}(s):= \suma_\alpha \tauec_{k*}(\alpha) N^{-s}(\alpha),
&\qquad&
M_{k*}(x):= \suma_{N(\alpha) \le x} \tauec_{k*}(\alpha).
\end{eqnarray*}

\begin{lemma}\label{th:tauec_dirichlet_series}
\begin{eqnarray}
\label{eq:dirichlet*-real}
G_{k*}(s) &=&
	\zeta(s)
	\zeta^{(k^2 + k - 2)/2}(2s)
	\zeta^{(-k^2 + k)/2}(3s)
	\zeta^{(-k^4 + 7k^2 - 6k)/12}(4s) \times
	\\ \nonumber
	&\times&
	\zeta^{(5k^4 - 6k^3 - 5k^2 + 6k)/24}(5s)
	g_{k*}(s), \\
\label{eq:dirichlet}
F_k(s)   &=&
	Z(s)
	Z^{k-1}(2s)
	Z^{(k-k^2)/2}(5s)
	Z^{(-k^3 + 6k^2 - 5k)/6}(6s)
	\times
	\\ \nonumber
	&\times&
	Z^{(k^3 - 4k^2 + 3k)/2}(7s)
	Z^{(3k^4 - 26k^3 + 57k^2 - 34k)/24}(8s)
	f_k(s), \\
\label{eq:dirichlet*}
F_{k*}(s) &=&
	Z(s)
	Z^{(k^2 + k - 2)/2}(2s)
	Z^{(-k^2 + k)/2}(3s)
	Z^{(-k^4 + 7k^2 - 6k)/12}(4s) \times
	\\ \nonumber
	&\times&
	Z^{(5k^4 - 6k^3 - 5k^2 + 6k)/24}(5s)
	f_{k*}(s),
\end{eqnarray}
where Dirichlet series $f_k(s)$ are absolutely convergent for $\Re s>1/9$ and Dirichlet series for $f_{k*}(s)$, $g_{k*}(s)$ are absolutely convergent for $\Re s>1/6$.
\end{lemma}

\proof
The statements can be directly verified with the help of the Bell series for corresponding functions. For example, for $\tauec_k$ we have following representation:
\begin{multline*}
\left( \sum_{a=0}^\infty \tauec_k(\p^a) x^a \right) (1-x)  (1-x^2)^{k-1}  (1-x^5)^{(k-k^2)/2}  (1-x^6)^{(-k^3 + 6k^2 - 5k)/6} \times
\\
\times (1-x^7)^{(k^3 - 4k^2 + 3k)/2}  (1-x^8)^{(3k^4 - 26k^3 + 57k^2 - 34k)/24} =  1+O(x^9) .
\end{multline*}
Then \eqref{eq:dirichlet} follows from the representation of $F_k$ and $Z$ in the form of infinite products by $\p$:
$$
F_k(s) = \prod_\p \left( \sum_{a=0}^\infty \tauec_k(\p^a) x^a \right),
\qquad
Z(s) = \prod_\p (1-\p^{-s})^{-1}.
$$
Identities \eqref{eq:dirichlet*-real} and \eqref{eq:dirichlet*} can be proved the same way.
\endproof

Let us define ${\mathbf a} := (1, \underbrace{2,\ldots,2}_{l})$,
$$
\tau({\mathbf a}; n) := \sum_{d_0d_1^2\cdots d_l^2=n} 1, \qquad T({\mathbf a}; x) := \sum_{n\le x} \tau({\mathbf a};n) = \sum_{d_0d_1^2\cdots d_l^2 \le x} 1,$$
Due to \cite[Th.~6.10]{kratzel1988} we have
\begin{equation}\label{eq:tau122-sum}
T({\mathbf a}; x) = C_1 x + x^{1/2} Q(\log x) + O(x^{w_l+\eps}),
\end{equation}
where $Q$ is a polynomial with computable coefficients, $\deg Q = l-1$, and $w_l \le (2l+1)/(4l+5)$. For some special values of $l$ better estimates of the error term can be obtained. For example, $w_1 \le 1057/4785$ (see \cite{graham1988}) and $w_2 \le 8/25$ due to \cite[(6.16)]{kratzel1988}.

\begin{theorem}\label{th:taue_sum}
\begin{equation*}\label{eq:taue_sum}
T_{k*}(x) = A_k x + x^{1/2} P_k(\log x) O(x^{v_k+\eps}),
\end{equation*}
where $P_k$ is a polynomial with computable coefficients, $\deg P_k = (k^2+k-4)/2$, and
$$ v_k = \max(w_{(k^2+k-2)/2}, 1/3). $$
\end{theorem}

\proof
Let $l=(k^2+k-2)/2$.
Identity \eqref{eq:dirichlet*-real} implies
\begin{equation}\label{eq:taue-convolution}
\taue_{k*} = \tau({\mathbf a};\cdot) \star f, \qquad T_{k*}(x) = \sum_{n\le x} T({\mathbf a}; x/n) f(n),
\end{equation}
where
series $\sum_{n=1}^\infty f(n) n^{-\sigma}$ are absolutely convergent for $\sigma>1/3$.

One can plainly estimate:
\begin{eqnarray}
\label{eq:f1-sum-tail}
\sum_{n>x} {f(n)\over n} &\ll& x^{-2/3+\eps} \sum_{n>x} {f(n)\over n^{1/3+\eps}} \ll x^{-2/3+\eps},
\\
\label{eq:f2-sum-tail}
\sum_{n>x} {f(n)  \log^a n \over n^{1/2}} &\ll& x^{-1/6+\eps} \sum_{n>x} {f(n) \log^a n\over n^{1/3+\eps}} \ll x^{-1/6+\eps}.
\end{eqnarray}
Substituting estimates \eqref{eq:tau122-sum}, \eqref{eq:f1-sum-tail} and \eqref{eq:f2-sum-tail} into \eqref{eq:taue-convolution} we get
\begin{multline*}
T_{k*}(x)
= C_1 x \sum_{n\le x} {f(n)\over n}
+ x^{1/2} \sum_{n\le x} {f(n) Q(\log(x/n))\over n^{1/2}}
+ O(x^{w_l+\eps}) + O(x^{1/3+\eps})
= \\
= A_k x + x^{1/2} P_k(\log x) + O(x^{v_k+\eps}).
\end{multline*}

\endproof

\begin{lemma}\label{th:tauec_res}
\begin{equation}\label{eq:tauec_res}
\res_{s=1} F_k(s)x^s/s = C_kx,
\qquad
\res_{s=1} F_{k*}(s)x^s/s = C_{k*}x,
\end{equation}
where
\begin{eqnarray}
\label{eq:C}
C_{k}   &= & \displaystyle {\pi\over 4} \prod_\p \left( 1 + \sum_{a=2}^\infty {\tau_k(a)-\tau_k(a-1) \over N^{a}(\p)} \right) ,
\\
\label{eq:C*}
C_{k*} &= & \displaystyle {\pi\over 4} \prod_\p \left( 1 + \sum_{a=2}^\infty {\tt_k(a)-\tt_k(a-1) \over N^{a}(\p)} \right).
\end{eqnarray}
\end{lemma}

\proof
As a consequence of the representation \eqref{eq:dirichlet} we have
$$ {F_k(s) \over Z(s)}
= \prod_p \left( 1 + \sum_{a=1}^\infty {\tau_k(a) \over N^{as}(\p)} \right) (1-\p^{-1})
= \prod_\p \left( 1 + \sum_{a=2}^\infty {\tau_k(a)-\tau_k(a-1) \over N^{as}(\p)} \right),
$$
and so function $F_k(s)/Z(s)$
is regular in the neighbourhood of $s=1$. At the same time we have
$$\res_{s=1} Z(s) = L(1, \chi_4) \res_{s=1} \zeta(s) = {\pi\over 4}, $$
which implies \eqref{eq:C}. The proof of \eqref{eq:C*} is similar.
\endproof

Numerical values of $C_k$ and $C_{k*}$ can be calculated in PARI/GP~\cite{parigp} with the use of the transformation
$$
\prod_\p f\bigl(N(\p)\bigr)
=
f(2)
\prod_{p=4k+1} f(p)^2
\prod_{p=4k+3} f(p^2)
$$
due to Lemma  \ref{l:gaussian-primes}. For example,
$$ C_2 \approx 1{,}156\,101, \qquad C_{2*}  \approx 1{,}524\,172. $$

\begin{theorem}\label{th:tauec_sum}
\begin{eqnarray}
\label{eq:M}
M_k(x) &=& C_kx + O(x^{1/2}\log^{3+4(k-1)/3} x),
\\
\label{eq:M*}
M_{k*}(x) &=& C_{k*}x + O(x^{1/2}\log^{3+2(k^2+k-2)/3} x),
\end{eqnarray}
where $C_k$ and $C_{k*}$ were defined in \eqref{eq:C} and \eqref{eq:C*}.
\end{theorem}

\proof
By Perron formula and by \eqref{eq:max-order-simple} for $c=1+1/\log x$, $\log T \asymp \log x$ we have
$$ M_k(x) = {1\over 2\pi i} \int_{c-iT}^{c+iT} F_k(s) {x^s\over s} ds + O\left( x^{1+\eps} \over T \right). $$
Suppose $d=1/2-1/\log x$. Let us shift the interval of integration to $[d-iT, d+iT]$. To do this consider an integral about a closed rectangle path with vertexes in
$$d-iT, ~ d+iT, ~ c+iT \text{~and~} c-iT.$$
There are two poles in $s=1$ and $s=1/2$ inside the contour. The residue at $s=1$ was calculated in  \eqref{eq:tauec_res}. The residue at $s=1/2$ is equal to $Dx^{1/2}$, $D=\const$ and will be absorbed by error term (see below).

Identity \eqref{eq:dirichlet} implies
$$F_k(s) = Z(s) Z^{k-1}(2s) f_k(s), $$
where $f_k(s)$ is regular for $\Re s > 1/3$, so for each $\eps>0$ it is uniformly bounded for~$\Re s > 1/3+\eps$.

Let us estimate the error term using Lemma \ref{l:zeta-moments} and Lemma \ref{l:phragmen}. The error term absorbs values of integrals about three sides of the integration's rectangle. We take into account $Z(s) = \zeta(s) L(s, \chi_4)$. On the horizontal segments we have
\begin{eqnarray*}
\int_{d+iT}^{c+iT} Z(s) Z^{k-1}(2s) {x^s\over s} ds
&\ll&
\max_{\sigma \in [d,c]} Z(\sigma+iT) Z^{k-1}(2\sigma + 2iT) x^\sigma T^{-1}
\ll \\ &\ll&
x^{1/2} T^{2\theta-1} \log^{4(k-1)/3}T + x T^{-1} \log^{4/3}T,
\end{eqnarray*}

It is well-known that $\zeta(s) \sim (s-1)^{-1}$ in the neighborhood of $s=1$.
So on~$[d, d+i]$ we get
\begin{multline*}
\int_{d}^{d+i} Z(s) Z^{k-1}(2s) {x^s\over s} ds
\ll
x^{1/2} \int_0^1 {\zeta^{k-1}(2d+2it)} dt
\ll \\
\ll
x^{1/2} \int_0^1 {dt\over |it-1/\log x|^{k-1}}
\ll
x^{1/2} \log^{k-1} x,
\end{multline*}
and for the rest of the vertical segment we have
\begin{multline*}
\int_{d+i}^{d+iT} Z(s) Z^{k-1}(2s) {x^s\over s} ds
\ll \\ \ll
\Bigl(
\int_1^T |\zeta(1/2+it)|^4 {dt\over t}
\int_1^T |L(1/2+it, \chi_4)|^4 {dt\over t}
\Bigr)^{1/4}
\Bigl(
\int_1^T |Z(1+2it)|^{2(k-1)} {dt\over t}
\Bigr)^{1/2}
\ll \\ \ll
x^{1/2} (\log^5 T \cdot \log^{8(k-1)/3+1} T)^{1/2}
\ll
x^{1/2} \log^{3+4(k-1)/3} T.
\end{multline*}
The choice $T=x^{1/2+\eps}$ finishes the proof of \eqref{eq:M}.

The proof of \eqref{eq:M*} is similar, but due to \eqref{eq:dirichlet*} one have replace $k-1$ by $(k^2+k-2)/2$.
\endproof

\bibliographystyle{ugost2008s}
\bibliography{taue}

\end{document}